\overfullrule=0pt
\centerline {\bf Weights sharing the same eigenvalue}\par
\bigskip
\bigskip
\centerline {BIAGIO RICCERI}\par
\bigskip
\bigskip
\centerline {\it Dedicated to Professor Wataru Takahashi, with esteem and friendship, on his seventieth birthday}\par
\bigskip
\bigskip
{\bf Abstract:} Here is the simplest particular case of our main result: let $f:{\bf R}\to {\bf R}$ be a function of class
$C^1$, with $\sup_{\bf R}f'>0$,  such that 
$$\lim_{|\xi|\to +\infty}{{f(\xi)}\over {\xi}}=0\ .$$
Then, for each $\lambda>{{\pi^2}\over {\sup_{\bf R}f'}}$, the set of all $u\in H^1_0(]0,1[)$ for which the problem
$$\cases{-v''=\lambda f'(u(x)) v & in $]0,1[$\cr & \cr 
v(0)=v(1)=0\cr}$$
has a non-zero solution is closed and not $\sigma$-compact in $H^1_0(]0,1[)$.\par
\bigskip
\bigskip
{\bf Key words:} Eigenvalue; weight; Dirichlet problem; $\sigma$-compact.\par
\bigskip
\bigskip
{\bf Mathematics Subject Classification:} 35P05; 47A75; 58C40.\par
\bigskip
\bigskip
\bigskip
\bigskip

{\bf 1. Introduction}\par
\bigskip
\bigskip
Let $\Omega\subset {\bf R}^n$ be a smooth bounded domain. We consider the Sobolev space $H^1_0(\Omega)$ endowed
with the scalar product
$$\langle u,v\rangle = \int_{\Omega}\nabla u(x)\nabla v(x)dx$$
and the induced norm
$$\|u\|=\left ( \int_{\Omega}|\nabla u(x)|^2dx\right ) ^{1\over 2}\ .$$
We are interested in pairs $(\lambda,\beta)$, where $\lambda$ is a positive number and $\beta$ is
a measurable function, such that the linear problem
$$\cases{-\Delta v=\lambda \beta(x) v & in $\Omega$\cr & \cr 
v_{|\partial \Omega}=0\cr}$$
has a non-zero weak solution, that is to say a $v\in H^1_0(\Omega)\setminus \{0\}$ such
that 
$$\int_{\Omega}\nabla v(x)\nabla w(x)dx=\lambda\int_{\Omega}\beta(x)v(x)w(x)dx$$
for all $w\in H^1_0(\Omega)$.\par
\smallskip
If this happens, we say that $\lambda$ is an eigenvalue related to the weight $\beta$.\par
\smallskip
While, the structure of the set of all eigenvalues related to a fixed weight is well understood,
it seems that much less is known about 
the structure of the set of all weights $\beta$ for which a fixed
positive number $\lambda$ turns out to be an eigenvalue related to $\beta$. \par
\smallskip
In this very short note, we intend to give a contribution along the latter direction. \par
\smallskip
More precisely,
we identify a quite general class of continuous functions $g: {\bf R}\to
{\bf R}$ such that, for each $\lambda$ in a suitable interval, the set of all $u\in H^1_0(\Omega)$
for which $\lambda$ is an eigenvalue related to the weight $g(u(\cdot))$ is closed and
not $\sigma$-compact in $H^1_0(\Omega)$.\par
\bigskip
\bigskip
{\bf 2. Results}\par
\bigskip
\bigskip
Let us recall that a set in a topological space is said to be $\sigma$-compact if it is
the union of an at most countable family of compact sets.\par
\smallskip
For each $\alpha\in L^{\infty}(\Omega)\setminus \{0\}$, with $\alpha\geq 0$, we denote by $\lambda_{\alpha}$ 
the first eigenvalue of the problem
$$\cases{-\Delta v=\lambda \alpha(x) v & in $\Omega$\cr & \cr 
v_{|\partial \Omega}=0\ .\cr}$$
Let us recall that
$$\lambda_{\alpha}=\min_{v\in H^1_0(\Omega)\setminus \{0\}}{{\|v\|^2}\over {\int_{\Omega}\alpha(x)|v(x)|^2dx}}\ .$$
With the conventions ${{1}\over {+\infty}}=0$, ${{1}\over {0}}=+\infty$, 
here is the statement of the result introduced above:\par
\medskip
THEOREM 1. - {\it Let $f:{\bf R}\to {\bf R}$ be a function of class $C^1$ such that
$$\max\left \{ 0, 2\limsup_{|\xi|\to +\infty}{{\int_0^{\xi}f(t)dt}\over {\xi^2}},
\limsup_{|\xi|\to +\infty}{{f(\xi)}\over {\xi}}\right \}<\sup_{{\bf R}}f'\ .$$
Moreover, if $n\geq 2$, assume that 
$$\sup_{\xi\in {\bf R}}{{|f'(\xi)|}\over {1+|\xi|^q}}<+\infty$$
for some $q>0$, with $q<{{4}\over {n-2}}$ if $n\geq 3$.\par
Then, for each $\alpha\in L^{\infty}(\Omega)\setminus \{0\}$, with $\alpha\geq 0$,
and for every $\lambda$ satisfying
$${{\lambda_{\alpha}}\over {\sup_{{\bf R}}f'}}<\lambda<{{\lambda_{\alpha}}\over {
\max\left \{ 0, 2\limsup_{|\xi|\to +\infty}{{\int_0^{\xi}f(t)dt}\over {\xi^2}},
\limsup_{|\xi|\to +\infty}{{f(\xi)}\over {\xi}}\right \}}}$$
the set of all $u\in H^1_0(\Omega)$ for which the problem
$$\cases{-\Delta v=\lambda \alpha(x)f'(u(x))v & in $\Omega$\cr & \cr 
v_{|\partial \Omega}=0\cr}$$
has a non-zero weak solution is closed and not $\sigma$-compact in $H^1_0(\Omega)$.}\par
\medskip
REMARK 1. - It is worth noticing that the linear hull of any closed and not $\sigma$-compact set in
$H^1_0(\Omega)$ is infinite-dimensional. This comes from the fact that any closed set in a finite-dimensional
normed space is $\sigma$-compact.\par
\medskip
The key tool we use to prove Theorem 1 is Theorem 2 below whose proof, in turn, is entirely
based on the following particular case of a result recently established in [1]:\par
\medskip
THEOREM A ([1], Theorem 2) {\it Let $(X,\langle\cdot,\cdot\rangle)$ be an infinite-dimensional real Hilbert space and let $I:X\to {\bf R}$
be a sequentially weakly lower semicontinuous, not convex functional of class $C^2$ such that
$I'$ is closed and $\lim_{\|x\|\to +\infty} (I(x)+\langle z,x\rangle) =+\infty$ for all $z\in X$.\par
Then, the set
$$\{x\in X : I''(x)\hskip 3pt is\hskip 3pt not\hskip 3pt invertible\}$$
is closed and not $\sigma$-compact.}\par
\medskip
THEOREM 2. - {\it Let $(X,\langle\cdot,\cdot\rangle)$ be an infinite-dimensional real Hilbert space,  
and let $J:X\to {\bf R}$ be a functional of class $C^2$, with compact derivative. 
For each $\lambda\in {\bf R}$,  put
$$A_{\lambda}=\{x\in X : y=\lambda J''(x)(y)\hskip 5pt for\hskip 5pt some\hskip 5pt y\in X\setminus \{0\}\}\ .$$
Assume that
$$\max\left\{0,2\limsup_{\|x\|\to +\infty}{{J(x)}\over {\|x\|^2}},
\limsup_{\|x\|\to +\infty}{{\langle J'(x),x\rangle}\over {\|x\|^2}}\right\}<\sup_{(x,y)\in X\times (X\setminus \{0\})}{{\langle J''(x)(y),y\rangle}\over
{\|y\|^2}}\ .$$
Then,  for every $\lambda$ satisfying
$$ \inf_{\{(x,y)\in X\times X : \langle J''(x)(y),y\rangle>0\}}{{\|y\|^2}\over {\langle J''(x)(y), y\rangle}}<\lambda<{{1}\over 
{\max\left\{0,2\limsup_{\|x\|\to +\infty}{{J(x)}\over {\|x\|^2}},
\limsup_{\|x\|\to +\infty}{{\langle J'(x),x\rangle}\over {\|x\|^2}}\right\}}}\ ,$$
 the set
$A_{\lambda}$ is closed and not $\sigma$-compact\ .}\par
\smallskip
PROOF. Fix $\lambda$ satisfying
$$ \inf_{\{(x,y)\in X\times X : \langle J''(x)(y),y\rangle>0\}}{{\|y\|^2}\over {\langle J''(x)(y), y\rangle}}<\lambda<{{1}\over 
{\max\left\{0,2\limsup_{\|x\|\to +\infty}{{J(x)}\over {\|x\|^2}},
\limsup_{\|x\|\to +\infty}{{\langle J'(x),x\rangle}\over {\|x\|^2}}\right\}}}\ .$$
 For each $x\in X$, put
$$I_{\lambda}(x)={{1}\over {2}}\|x\|^2-\lambda J(x)\ .$$
Clearly, for some $(x,y)\in X\times X$, with $\langle J''(x)(y),y\rangle>0$,  we have
$$\left \langle y-\lambda J''(x)(y), y\right \rangle <0$$
and so, since 
$$I''_{\lambda}(x)(y)= y-\lambda J''(x)(y)\ ,$$
by a classical characterization (Theorem 2.1.11 of [2]), the functional $I_{\lambda}$ is not convex.
Now, let us show that
$$\lim_{\|x\|\to +\infty}\|x-\lambda J'(x)\|=+\infty\ .\eqno{(1)}$$
Indeed, for each $x\in X\setminus \{0\}$, we have
$$\|x-\lambda J'(x)\|=\sup_{\|y\|=1}\langle x+\lambda J'(x),y\rangle\geq \left\langle x-\lambda J'(x),{{x}\over {\|x\|}}
\right\rangle\geq
\|x\|\left ( 1-\lambda{{\langle J'(x),x\rangle}\over {\|x\|^2}}\right )\ .\eqno{(2)}$$
On the other hand, we also have
$$\liminf_{\|x\|\to +\infty}\left ( 1-\lambda{{\langle J'(x),x\rangle}\over {\|x\|^2}}\right )=
1-\lambda\limsup_{\|x\|\to +\infty}{{\langle J'(x),x\rangle}\over {\|x\|^2}}>0\ .\eqno{(3)}$$
So, $(1)$ is a direct consequence of $(2)$ and $(3)$.
Furthermore,  for each $z\in X$, since
$$I_{\lambda}(x)+\langle z,x\rangle=\|x\|^2
\left ( {{1}\over {2}}-\lambda{{J(x)}\over {\|x\|^2}}+{{\langle z,x\rangle}\over {\|x\|^2}}\right )$$
and
$$\liminf_{\|x\|\to +\infty}\left ( {{1}\over {2}}-\lambda{{J(x)}\over {\|x\|^2}}+{{\langle z,x\rangle}\over {\|x\|^2}}\right )=
{{1}\over {2}}-\lambda\limsup_{\|x\|\to +\infty}{{J(x)}\over {\|x\|^2}}>0\ ,$$
we have
$$\lim_{\|x\|\to +\infty}(I_{\lambda}(x)+\langle z,x\rangle)=+\infty\ .$$
Since $J'$ is compact, on the one hand, $J$ is sequentially
weakly continuous ([4], Corollary 41.9) and, on the other hand, in view of $(1)$, the operator
$I_{\lambda}'$ is closed ([3], Example 4.43).
 The compactness of $J'$ also implies that, for each
$x\in X$, the operator $J''(x)$ is compact ([3], Proposition 7.33) and so, for each $\lambda\in {\bf R}$,
 the operator $y\to y-\lambda J''(x)(y)$ is injective if and only if it is surjective ([3], Example 8.16).
 At this point, the fact that $A_{\lambda}$ is closed and not $\sigma$-compact follows directly from
Theorem A which can be applied to the functional $I_{\lambda}$.\hfill $\bigtriangleup$
\par
\medskip
{\it Proof of Theorem 1.}  Fix $\alpha\in L^{\infty}(\Omega)\setminus \{0\}$, with $\alpha\geq 0$.
For each $u\in H^1_0(\Omega)$, put
$$J_{f}(u)=\int_{\Omega}\alpha(x)F(u(x))dx\ ,$$
where
$$F(\xi)=\int_{0}^{\xi}f(t)dt\ .$$
Our assumptions ensure that the functional $J_f$ is of class $C^2$ in $H^1_0(\Omega)$, and we have 
$$\langle J_f'(u),v\rangle=\int_{\Omega}\alpha(x)f(u(x))v(x)dx\ ,$$
$$\langle J_f''(u)(v),w\rangle=\int_{\Omega}\alpha(x)f'(u(x))v(x)w(x)dx$$
for all $u,v,w\in H^1_0(\Omega)$. Moreover, $J_f'$ is compact. Fix $\nu>\limsup_{|\xi|\to +\infty}{{F(\xi)}\over
{\xi^2}}$. Then, for a suitable constant $c>0$, we have
$$F(\xi)\leq \nu\xi^2+c$$
for all $\xi\in {\bf R}$. Hence, for each $u\in H^1_0(\Omega)$, we obtain
$$J_f(u)\leq \nu\int_{\Omega}\alpha(x)|u(x)|^2dx+c\int_{\Omega}\alpha(x)dx\leq
\nu\lambda_{\alpha}^{-1}\|u\|^2+ c\int_{\Omega}\alpha(x)dx\ .$$
This clearly implies that
$$\limsup_{\|u\|\to +\infty}{{J_f(u)}\over {\|u\|^2}}\leq \lambda_{\alpha}^{-1}\limsup_{|\xi|\to +\infty}{{F(\xi)}\over
{\xi^2}}\ .\eqno{(4)}$$
In the same way, we obtain
$$\limsup_{\|u\|\to +\infty}{{\langle J_f'(u),u\rangle}\over {\|u\|^2}}\leq 
\lambda_{\alpha}^{-1}\limsup_{|\xi|\to +\infty}{{f(\xi)}\over
{\xi}}\ .\eqno{(5)}$$
Now, fix a function $\tilde v\in H^1_0(\Omega)$, with $\|\tilde v\|=1$,  such that
$$\int_{\Omega}\alpha(x)|\tilde v(x)|^2dx=\lambda_{\alpha}^{-1}\ .$$
Fix also $\epsilon>0$, $\tilde\xi\in {\bf R}$, with $f'(\tilde\xi)>0$, and a closed set
$C\subseteq \Omega$ so that
$$\int_C\alpha(x)|\tilde v(x)|^2dx>\int_{\Omega}\alpha(x)|\tilde v(x)|^2dx-\epsilon$$
and
$$\int_{\Omega\setminus C}\alpha(x)|\tilde v(x)|^2dx<{{\epsilon}\over {\sup_{[-|\tilde\xi|,|\tilde\xi|]}|f'|}}\ .$$
Finally, fix a function $\tilde u\in H^1_0(\Omega)$ such that
$$\tilde u(x)=\tilde\xi$$
for all $\xi\in C$ and
$$|\tilde u(x)|\leq |\tilde \xi|$$
for all $\xi\in \Omega$. Then, we have
$$f'(\tilde\xi)\left(\int_{\Omega}\alpha(x)|\tilde v(x)|^2dx-\epsilon\right) < f'(\tilde\xi)\int_{C}\alpha(x)|\tilde v(x)|^2dx=
\int_{\Omega}\alpha(x)f'(\tilde u(x))|\tilde v(x)|^2dx-\int_{\Omega\setminus C}\alpha(x)f'(\tilde u(x))|\tilde v(x)|^2dx$$
$$\leq \sup_{(u,v)\in H^1_0(\Omega)\times (H^1_0(\Omega)\setminus \{0\})}{{\int_{\Omega}\alpha(x)f'(u(x))|v(x)|^2dx}\over
{\|v\|^2}}+\epsilon\ .$$
Since $\tilde\xi$ and $\epsilon$ are arbitrary, we then infer that
$$\lambda_{\alpha}^{-1}\sup_{\bf R}f'\leq 
\sup_{(u,v)\in H^1_0(\Omega)\times (H^1_0(\Omega)\setminus \{0\})}{{\int_{\Omega}\alpha(x)f'(u(x))|v(x)|^2dx}\over
{\|v\|^2}}\ .\eqno{(6)}$$
Consequently, putting $(4)$, $(5)$ and $(6)$ together, we obtain
$$\max\left\{0,2\limsup_{\|u\|\to +\infty}{{J_f(u)}\over {\|u\|^2}},{{\langle J_f'(u),u\rangle}\over {\|u\|^2}}\right\}
\leq \lambda_{\alpha}^{-1}\max\left \{ 0, 2\limsup_{|\xi|\to +\infty}{{\int_0^{\xi}f(t)dt}\over {\xi^2}},
\limsup_{|\xi|\to +\infty}{{f(\xi)}\over {\xi}}\right \}$$
$$<\lambda_{\alpha}^{-1}\sup_{{\bf R}}f'\leq
\sup_{(u,v)\in H^1_0(\Omega)\times (H^1_0(\Omega)\setminus \{0\})}{{\int_{\Omega}\alpha(x)f'(u(x))|v(x)|^2dx}\over
{\|v\|^2}}\ .$$
Therefore, we can apply Theorem 2 taking $X=H^1_0(\Omega)$ and
$J=J_f$. Therefore, for every $\lambda$ satisfying 
$${{\lambda_{\alpha}}\over {\sup_{{\bf R}}f'}}<\lambda<{{\lambda_{\alpha}}\over {
\max\left \{ 0, 2\limsup_{|\xi|\to +\infty}{{\int_0^{\xi}f(t)dt}\over {\xi^2}},
\limsup_{|\xi|\to +\infty}{{f(\xi}\over {\xi}}\right \}}}\ ,$$ the
set $A_{\lambda}$ (defined in Theorem 2) is closed and not $\sigma$-compact in
$H^1_0(\Omega)$. But, clearly, a $u\in H^1_0(\Omega)$ belongs to $A_{\lambda}$
if and only if the problem
$$\cases{-\Delta v=\lambda \alpha(x)f'(u(x))v & in $\Omega$\cr & \cr 
v_{|\partial \Omega}=0\cr}$$
has a non-zero weak solution, and the proof is complete.\hfill $\bigtriangleup$
\bigskip
\bigskip
\noindent
{\bf Acknowledgement.} The author has been supported by the Gruppo Nazionale per 
l'Analisi Matematica, la Probabilit\`a e le loro Applicazioni (GNAMPA) of the Istituto Nazionale di Alta Matematica (INdAM).

\vfill\eject
\centerline {\bf References}\par
\bigskip
\bigskip
\noindent
[1]\hskip 5pt B. RICCERI, {\it Singular points of non-monotone potential operators}, preprint.\par
\smallskip
\noindent
[2]\hskip 5pt C. Z\u{A}LINESCU, {\it Convex analysis in general vector spaces}, World Scientific, 2002.\par
\smallskip
\noindent
[3]\hskip 5pt E. ZEIDLER, {\it Nonlinear functional analysis and its
applications}, vol. I, Springer-Verlag, 1986.\par
\smallskip
\noindent
[4]\hskip 5pt E. ZEIDLER, {\it Nonlinear functional analysis and its
applications}, vol. III, Springer-Verlag, 1985.
\bigskip
\bigskip
\bigskip
\bigskip
Department of Mathematics\par
University of Catania\par
Viale A. Doria 6\par
95125 Catania\par
Italy\par
{\it e-mail address}: ricceri@dmi.unict.it

\bye
Set
$$f(x,\xi)=\beta(x)\psi(\xi)$$
for all $(x,\xi)\in \Omega\times {\bf R}$. Clearly, $f\in {\cal A}$.
Let us apply Theorem 3, taking $X=H^1_0(\Omega)$ and
$J=-I_f$. Our assumptions imply that $J$ is of class $C^2$, with
$$\langle J''(u)(v),w\rangle=-\int_{\Omega}\psi'(u(x))v(x)w(x)dx$$
for all $u,v,w\in X$. Thanks to $(20)$, from the proof of Proposition 3,
we know that
$$\lim_{\|u\|\to +\infty}{{\|I_{f}'(u)\|}\over {\|u\|}}=0 \ .\eqno{(21)}$$
We claim that
$$\lim_{\|u\|\to +\infty}{{I_f(u)}\over {\|u\|^2}}=0 \ .$$
Arguing by contradiction, assume that there exist $\gamma>0$ and a sequence
  $\{u_k\}$, with $\lim_{k\to \infty}\|u_k\|=+\infty$, such that
$$I_f(u_k)>\gamma \|u_k\|^2$$
for all $k\in {\bf N}$. By the mean value theorem, there is $t_k\in [0,1]$
so that $I_f(u_k)=\langle I_{f}'(t_k u_k),u_k\rangle$. So, we clearly
have
$$\|I_{f}'(t_k u_k)\|>\gamma\|u_k\|\ .\eqno{(22)}$$
Since $I_{f}'$ is bounded on each bounded subset of $X$, we then infer that,
up to a subsequence, $\lim_{k\to \infty}\|t_k u_k\|=+\infty$. Consequently,
if $k$ is large enough, in view of $(21)$, we would have
$$\|I_{f}'(t_k u_k)\|<\gamma\|t_k u_k\|$$
which contradicts $(22)$. So, in view of Propositions 2 and 3, all the
assumptions of Theorem 3 are satisfied. The conclusion then follows directly
from Theorem 3, taking Remark 2 into account.\hfill $\bigtriangleup$\par

\bye